\documentclass[notitlepage,11pt]{article}  

\usepackage{amssymb,amsmath,latexsym,epsfig,euscript}

\newlength{\jmr}
\newlength{\hwl}
\newlength{\khov}
\newtheorem{smale}{The Shub-Smale $\pmb{\tau}$-Conjecture}

\newtheorem{amd}{Arithmetic Multivariate Descartes' Rule (Special Case)}

\newtheorem{pcfew}{$p$-adic Complex Fewnomial Theorem}

\newtheorem{padic}{The $p$-adic Digit Conjecture}

\newtheorem{thm}{Theorem}
\newtheorem{prop}{Proposition}
\newtheorem{dfn}{Definition}

\newtheorem{rem}{Remark}

\newtheorem{ex}{Example}
\newtheorem{lemma}{Lemma}

\newcommand{\ord}{\mathrm{ord}}

\newcommand{\am}{{\mathbf{AM}}}

\newcommand{\np}{{\mathbf{NP}}}

\newcommand{\bpp}{{\mathbf{BPP}}}

\newcommand{\pp}{\mathbf{P}}

\newcommand{\thth}{{\underline{\mathrm{th}}}}

\newcommand{\nd}{{\underline{\mathrm{nd}}}}

\newcommand{\Q}{\mathbb{Q}}
\newcommand{\R}{\mathbb{R}}
\newcommand{\C}{\mathbb{C}}

\newcommand{\Z}{\mathbb{Z}}

\newcommand{\cp}{\mathfrak{p}}
\newcommand{\newt}{\mathrm{Newt}}

\newcommand{\conv}{\mathrm{Conv}}

\newcommand{\qed}{$\blacksquare$}
\newcommand{\dia}{$\diamond$}

\newcommand{\cL}{{\mathcal{L}}}

\newcommand{\bO}{\mathbf{O}}

\begin{document}
%

\title{\mbox{}\\
\vspace{-.8in}A Direct Ultrametric Approach to Additive Complexity and the 
Shub-Smale Tau Conjecture 
}
%

%

%
\author{J.\ Maurice Rojas\thanks{Partially supported by
NSF Grant DMS-0211458. }\\
Texas A\&M University\\
Department of Mathematics\\
TAMU 3368\\
College Station, TX \ 77843-3368\\
USA\\
e-mail: {\tt rojas@math.tamu.edu}\\
Web Page: {\tt  http://www.math.tamu.edu/\~{}rojas}\\
FAX: (979) 845-6028 }

\date{\today} 
\maketitle
\begin{abstract}
The Shub-Smale Tau Conjecture is a hypothesis relating 
the number of integral roots of a polynomial $f$ in one variable and 
the {\bf Straight-Line Program (SLP) complexity} of $f$. 
A consequence of the truth of this conjecture is that, 
for the Blum-Shub-Smale model over the complex numbers,  
P differs from NP. We prove two weak versions of the Tau Conjecture 
and in so doing show that the Tau Conjecture follows from an even 
more plausible hypothesis. 

Our results follow from a 
new $p$-adic analogue of earlier work relating 
real algebraic geometry to additive complexity. For instance, we can 
show that a nonzero univariate polynomial of {\bf additive complexity} $s$  
can have no more than $1+s^3(s+1)(7.5)^s s!=O(e^{s\log s})$ roots in the 
$2$-adic rational numbers $\Q_2$, thus dramatically improving an earlier 
result of the author. This immediately implies the same bound on the number 
of ordinary rational roots, whereas the\\ 
\scalebox{.85}[1]{best previous upper bound via 
earlier techniques from real algebraic geometry was a quantity in 
$\Omega((22.6)^{s^2}$).} 

This paper presents another step in the author's program 
of establishing an algorithmic arithmetic version of 
fewnomial theory. 
\end{abstract}




\section{Introduction}
We show how ultrametric geometry can be used to give new sharper 
lower bounds for the additive complexity of polynomials. As a 
consequence, we advance further a number-theoretic approach toward 
an analogue of the $\pp\!\neq\!\np$ conjecture for the BSS 
model over $\C$. Let us now clarify all these notions. 
(Simple explicit examples of our constructions appear in 
Section \ref{sub:related} below.) 

The {\bf Blum-Shub-Smale (BSS) model over $\C$} can be thought of 
roughly as a Turing machine augmented with additional registers that 
perform arithmetic operations on complex numbers at unit cost. 
(See \cite{bcss} for a much more complete description.) This model 
of computation possesses a natural analogue of the classical 
$\pp\!\stackrel{?}{=}\!\np$ question: the $\pp_\C\!\stackrel{?}{=}\!\np_\C$ 
question. This new analogue was introduced with the hope of (a) enabling 
researchers to enlarge their usual bag of tricks with the 
arsenal of number theory and analysis and (b) gaining new information 
on the original $\pp\!\stackrel{?}{=}\!\np$ question. 

Indeed, it is now known (see, e.g., \cite{shub,duke,bcss}) that 
$\pp_\C\!=\!\np_\C 
\Longrightarrow \np\!\subseteq\!\bpp$. Since the latter containment is 
widely disbelieved, this makes it plausible that 
$\pp_\C\!\neq\!\np_\C$. The implications of $\pp_\C\!\neq\!\np_\C$ 
for classical complexity are still not clear. However, via an 
examination of the Boolean parts of $\pp_\C$ and $\np_\C$ one can see 
that $\pp_\C\!\neq\!\np_\C$ would provide some evidence that 
$\bpp\!\neq\!\am$ \cite{koiran,dzh}. Better still, Mike Shub and 
Steve Smale were able 
to obtain an elegant number-theoretic statement that {\bf implies}  
$\pp_\C\!\neq\!\np_\C$ (see \cite{duke} or \cite[Thm.\ 3, Pg.\ 127]{bcss}). 
Recall that $\Z$ denotes the integers. 
\begin{dfn} 
Suppose a polynomial in one variable $f\!\in\!\Z[x_1]$ is expressed as
a sequence of the form $(1,x_1,f_2,\ldots,f_N)$, where $f_N\!=\!f(x_1)$,
$f_0\!:=\!1$, $f_1\!:=\!x_1$, and for all $i\!\geq\!2$ we have that $f_i$
is a sum, difference, or product of some pair of elements $(f_j,f_k)$
with $j,k\!<\!i$. Such a computational sequence is a {\bf              
straight-line program (SLP)}. We then let $\pmb{\tau(f)}$ denote the
smallest possible value of $N-1$, i.e., the smallest length for such a
computation of $f$. \dia
\end{dfn} 
\begin{smale}
There is an absolute constant $\kappa\!\geq\!1$ such that
for all\\ $f\!\in\!\Z[x_1]\setminus\{0\}$, the
number of distinct roots of $f$ in $\Z$ is no more than 
$(\tau(f)+1)^\kappa$. 
\end{smale}
It is easily checked that $\deg f\!\leq\!2^{\tau(f)}$ so we at 
least know that $f$ has no more than $2^{\tau(f)}$ integral roots.
However, little else is known about the relation of $\tau(f)$ to the 
number of integral roots of $f$:  
as of mid-2003, the $\tau$-conjecture remains open, even in the special 
case $\kappa\!=\!1$. The polynomial $(x-2^1)(x-2^2)\cdots (x-2^{2^j})$ 
easily shows that the $\tau$-conjecture fails if we allow $\kappa\!<\!1$. 

To state our first main result, let us introduce a family of weaker 
hypotheses depending on an integer parameter $p$.  
\begin{padic} 
There is an absolute constant $c_p\!>\!0$ such that 
for all $f\!\in\!\Z[x_1]\setminus\{0\}$, the
number of distinct roots $x\!\in\!\Z$ of $f$ with 
$x\!\equiv\!1 \ (\mathrm{mod} \ p)$ is no more than $(\tau(f)+1)^{c_p}$. 
\end{padic} 

Clearly, the $\tau$-conjecture implies the $p$-adic Digit 
Conjecture for any integer $p$. However, a strong converse holds as well. 
\begin{thm}
\label{thm:padic} 
Suppose the $p$-adic Digit Conjecture is true for some fixed prime $p$. 
Then $\pp_\C\!\neq\!\np_\C$. In particular, if the $p$-adic Digit Conjecture 
is true for some fixed prime $p$, then the $\tau$-conjecture is true. 
\end{thm} 
The last implication follows easily from our next main theorem. 
\begin{dfn}
\label{dfn:add} 
Letting $R$ be any ring, we say that 
$f\!\in\!R[x_1]$ has {\bf additive complexity $\pmb{\leq s}$ (over 
$R$)} iff there exist constants $c_1,d_1, 
\ldots,c_s$, $d_s,c_{s+1}\!\in\!R$ and arrays of
nonnegative integers $[m_{i,j}]$ and $\left[m'_{i,j}\right]$ 
with $f(x)\!=\!c_{s+1}\prod\limits^s_{i=0} X^{m_{i,s+1}}_i$,
where $X_0\!=\!x$, $X_1 = c_1X^{m_{0,1}}_0+
d_1X^{m'_{0,1}}_0$, 
and $X_j = c_j \left(\prod\limits^{j-1}_{i=0}
X^{m_{i,j}}_i\right) + d_j \left(\prod\limits^{j-1}_{i=0} X^{m'_{i,j}}_i
\right)$ for all $j\!\in\!\{2,\ldots,s\}$. 
We then define the {\bf additive complexity (over $\pmb{R}$) of
$\pmb{f}$}, $\pmb{\sigma_R(f)}$,
to be the least $s$ in such a presentation of $f$ as an algebraic
expression. Finally, additive complexity without appellation, and 
$\sigma(f)$ will be understood to mean $\sigma_\Z(f)$. \dia
\end{dfn}
Note in particular that repeated additions or subtractions in different 
sub-expressions are thus {\bf not} counted, e.g., 
$\sigma(9(x-7)^{99}(2x+1)^{43}-11(x-7)^{999}
(2x+1)^3)\!\leq\!3$. 
Recall that for any prime number $p\!\in\!\Z$, the symbols $\Z_p$, $\Q_p$, 
and $\C_p$ respectively denote the $p$-adic integers, the $p$-adic rational 
numbers, and the $p$-adic complex numbers (cf.\ Section \ref{sec:back}). 
\begin{dfn}
\label{dfn:padic}
Given any $p$-adic integer $x\!\in\!\Z_p$ define its {\bf $\pmb{p}$-adic 
valuation}, $\ord_p x$, to be the largest $s$ such that $p^s$ divides $x$ 
(so $\ord_p 0\!:=\!+\infty$). 
We then extend $\ord_p$ to $\Q_p$ by setting $\ord_p \frac{m}{n}\!:=\!\ord_p 
m - \ord_p n$, for any nonzero $m,n\!\in\!\Z_p$. Finally, we define the
{\bf $\pmb{p}$-adic norm} of any $p$-adic rational number $x$, $|x|_p$,
to be $p^{-\ord_p x}$. \dia 
\end{dfn}
\begin{thm}
\label{thm:abs} 
Let $p$ be any prime number and let $\#$ denote the operation of taking 
set cardinality. Also, following the notation above, 
let $N_p(s):=\max\limits_{\substack{f\in\Z[x_1]\\ \sigma(f)\leq s}}  
\#\left\{|x|_p \; | \; x\!\in\!\C_p\!\setminus\!\{0\} \text{ is a root 
of } f\right\}$. Then 
\[s\!\leq\!N_p(s)\!\leq\!s(s+1)/2.\] Furthermore, 
the lower bound can be realized even if one restricts to 
roots in $\Z$. In particular, the number of possible locations for 
the first nonzero base-$p$ digit of an integral root of $f$  
is no more than $s(s+1)/2$. 
\end{thm} 
\begin{rem} 
Alternatively, the quantity $N_p(s)$ can be thought of as the $p$-adic 
analogue of the maximum number of distinct absolute values 
that can occur for the roots of a fixed polynomial of bounded additive 
complexity. Note in particular that $(x+1)^d-1$ has additive 
complexity $\leq\!2$ but exactly $\left\lceil\frac{d}{2}\right\rceil$ 
distinct absolute values for its complex roots. So, differing radically from 
the analogous situation over $\C$, the maximum number of distinct 
absolute values over $\C_p$ is not only a 
{\bf finite} function of $s$ and $p$, but sub-quadratic and independent of 
$p$. \dia 
\end{rem} 
Since $\Z\!\subset\!\Q\!\subset\!\Q_p\!\subset\!\C_p$, we can thus 
think of Theorem \ref{thm:abs} as a weakening of the $\tau$-conjecture where 
we count $p$-adic absolute values of roots instead of roots. The author 
believes that a lower bound quadratic in $s$ is also possible. 

Theorem \ref{thm:abs} in turn implies the truth of another weak version of 
the $\tau$-conjecture: Theorem \ref{thm:add} below, which may also be of 
independent interest in arithmetic circuit complexity. 
\begin{thm} 
\label{thm:add} 
Any polynomial $f\!\in\!\Z[x_1]\!\setminus\!\{0\}$ of 
additive complexity $\leq\!s$   
has no more than $O\left(pe^{s\log s}\right)$ 
roots in $\Q_p$. More precisely, $1+(p-1)s^2 
(7.5)^s s!N_p(s)$ suffices as an upper bound. In particular, the number of 
rational roots of an $f\!\in\!\Z[x_1]\!\setminus\!\{0\}$ with 
additive complexity $\leq\!s$ over $\Z$ is no more than 
$1+s^3(s+1)(7.5)^s s!$. 
\end{thm} 
\begin{rem}
Note that $\tau(f)$ is already unbounded for polynomials of the 
form $a+bx^d$, which clearly have additive complexity $1$. 
So $\sigma(f)$ can be tremendously smaller than $\tau(f)$. 
Our theorem above can thus be interpreted as a weak, but 
definitely non-trivial version of the $\tau$-conjecture. \dia 
\end{rem} 
Theorem \ref{thm:add} can be sharpened further (see 
Section \ref{sec:add} for the details) and even the coarse bound on 
rational roots above is the best to date in terms of additive complexity. 

Theorem \ref{thm:add} follows easily from a weakening of the $p$-adic 
Digit Conjecture that we actually can prove: Theorem \ref{thm:cx} 
below gives an upper bound on the number of roots of $f$ 
with ``first $p$-adic digit $1$'' which, while not polynomial in $\tau(f)$, 
is exponential in {\bf $\pmb{\sigma(f)}$, explicit, and applies to 
$\pmb{p}$-adic complex roots as well}. 
\begin{thm}
\label{thm:cx} 
For any $r\!>\!0$, let $C'_p(s,r)$ denote the maximum, 
over all $f\!\in\!\Z[x_1]$ with $\sigma(f)\!\leq \!s$, 
of \[\#\left\{\text{Roots of } f  \text{ in the closed complex
$p$-adic disk of radius } \frac{1}{p^r} \text{ about } 1 \right\}.\] 
Then $C'_p(s,r)\!\leq\!s^2s!\left(3+\frac{3}{r}\log_p\left(\frac{2}
{r\log p}\right)\right)^s$. 
\end{thm} 
Theorem \ref{thm:cx} can also be sharpened even further and the 
details are in Section \ref{sec:cx}.  The author suspects that 
the bound above can even be improved to 
$\left(1-\frac{\log r}{r}\right)^{O(s)}$.  

The proofs of Theorems \ref{thm:padic}, \ref{thm:abs}, \ref{thm:add}, 
and \ref{thm:cx} appear respectively in Sections \ref{sec:padic}, 
\ref{sec:abs}, \ref{sec:add}, and \ref{sec:cx}.  
In closing this introduction we note that all our preceding results 
admit natural extensions to $\cp$-adic fields (i.e., 
finite algebraic extensions of the usual $p$-adic numbers), 
roots of bounded degree over number fields (finite algebraic extensions 
of the ordinary rational nubmers), and even systems of multivariate 
polynomial equations. These extensions will be detailed in 
the full version of this paper. 

Let us now further detail earlier work on these problems as 
well as some background on $p$-adic numbers. 

\subsection{Related Work and More Examples} 
\label{sub:related}
Let us first point out that there are examples showing that 
we {\bf must} avoid $\R$ if we are to use field arithmetic 
tricks to solve the $\tau$-conjecture. 
\begin{ex} 
Consider the recurrence $g_{j+1}\!:=\!4g_j(1-g_j)$ with $g_1\!:=\!4x(1-x)$. 
It is then easily checked\footnote{ This example is well-known in dynamical 
systems, and the author thanks Gregorio Malajovich for pointing it out. }
that $g_j(x)-x$ has exactly $2^j$ roots in the open interval $(0,1)$, but 
$\tau(g_j(x)-x)\!=\!O(j)$. \dia 
\end{ex} 
The existence of an analogous example over the $p$-adic rationals 
is still an open question. 

That the $\tau$-conjecture is still open is a testament to the 
fact that we know far less about the complexity measures $\sigma$ 
and $\tau$ than we should. For example, there is still no 
more elegant method known to compute $\tau$ for a fixed 
polynomial than brute force 
enumeration. Also, the computability of additive complexity 
is still an open question, although a more efficient variant 
(allowing radicals as well) can be computed in triply exponential 
time \cite{grigokarp}. 

As for earlier earlier approaches to the $\tau$-conjecture, 
some investigations into $\tau$ were initiated by de Melo and 
Svaiter \cite{fux} (in the special case of constant polynomials) 
and Moreira \cite{gugu}. Unfortunately, not much more is known than 
(a) $\tau(f)$ is ``usually'' bounded below by $h(f)/\log h(f)$ 
(where $h(f)$ is $(1+\deg f)\max \log |c|$ and the maximum ranges 
over all coefficients $c$ of $f$) and (b) $\tau\!\leq\!1.1h(f)/\log h(f)$ 
for $h(f)$ sufficiently large \cite[Thm.\ 3]{gugu}. 

So information relating the {\bf integral} roots of $f$ with $\tau(f)$ 
was essentially 
non-existent, at least until independent work of Dima Yu.\ Grigoriev 
\cite{grigo} and Jean-Jacques Risler \cite{risler} that related additive 
complexity with the number of {\bf real} roots. (This was also preceded by 
important seminal work of Allan Borodin and Stephen A.\ Cook \cite{bocook} 
first proving the surprising fact that one could indeed bound 
the number of real roots in terms 
of additive complexity.) In particular, they showed that an 
$f\!\in\!\R[x_1]$ with $\sigma_\R(f)\!\leq\!s$ could have no more than 
$(s+2)^{3s+1}2^{\left.\left(9s^2 +5s+2\right)\right/2}$ 
(or $\Omega\left((22.6)^{s^2}\right)$) real roots 
\cite[Pg.\ 181, Line 6]{risler}. 

It was then discovered around 2001 that one could derive even sharper 
bounds by working $p$-adically. In particular, \cite{add} gives an upper 
bound of $1+s^2(22.5)^s s!$ for the number of $2$-adic rational roots 
of an $f\!\in\!\Q_2[x_1]$ with $\sigma_{\Q_2}(f)\!\leq\!s$, and 
states explicit generalizations to $\Q_p$, algebraic extensions, 
roots of bounded degree over number fields, 
and systems of multivariate polynomial equations. Our results thus 
improve the last bound by a factor close to $3^s$ and isolates the portion 
where the (necessary?) exponentiality is coming from. 

Historically, the results of \cite{add} were derived as a consequence of a 
higher-dimensional result, also due to the author, relying on a less 
efficient encoding: monomial term expansions. For completeness, 
we paraphrase the most relevant special case of the result below. 
\begin{amd}
\cite[Cor.\ 1 of Sec.\ 3]{amd}\\ Let $p$
be any prime and let $\Q^*_p\!:=\!\Q_p\setminus\{0\}$.
Suppose $\pmb{f_1},\ldots,\pmb{f_k}\!\in\!\Q_p[x^{\pm 1}_1,\ldots,x^{\pm 1}_n]
\setminus\{0\}$, $F\!:=\!(f_1,\ldots,f_k)$, $\pmb{m_i}$ is the total
number of distinct exponent
vectors appearing in $f_i$ (assuming all polynomials are written as
sums of monomials), and $N_i$ is the number of variables occuring in $f_i$.
Finally, let $m\!:=\!(m_1,\ldots,m_n)$, $N\!:=\!(N_1,\ldots,N_n)$,
and define $\pmb{B(\cL,m,N)}$ to be
the maximum number of isolated roots in $(\Q^*_p)^n$ of such an
$F$, counting multiplicities.\footnote{ For the sake of simplicity,
one can safely assume that the multiplicity of any isolated root is always
a positive integer. }
Then\\ $\pmb{B(\cL,m,N)\!\leq\!\left(\prod\limits^n_{i=1} 
(p-1)m_i(m_i-1)/2\right) \left\lfloor\prod\limits^n_{i=1} \left
\{c(m_i-1)N_i\left[1+\log_p
\left(\frac{(m_i-1)}{\log p} \right)\right]
\right\}\right\rfloor}$,
where $c\!:=\!\frac{e}{e-1}\!\leq\!1.582$. \qed
\end{amd}
This result generalized earlier work of Lenstra for the 
case of sparse polynomials in one variable \cite{lenstra2}, and 
was derived via an analogous result counting $p$-adic 
complex roots close to the point $(1,\ldots,1)$. 

The latter result, which provided the first effective version 
of earlier model-theoretic 
work of Denef, Lipshitz, and van den Dries \cite{vandenef,lipshitz}, 
can also be viewed as an arithmetic extension \cite{amd} of
Khovanski's Theorem on Complex Fewnomials \cite{few}.
\begin{pcfew}
\cite[Thm.\ 2 of Sec.\ 1.1]{amd}\\ Following the notation above,
suppose instead now that $k\!=\!n$ and that the coefficients of $F$ can lie
in $\C_p$. Let $r_1,\ldots,r_n\!>\!0$, $r\!:=\!(r_1,\ldots,r_n)$, and
let $C_p(m,N,r)$ denote the maximum number of geometrically isolated roots
$(x_1,\ldots,x_n)\!\in\!\C^n_p$ of $F$ with $\ord_p(x_i-1)\!\geq\!r_i$
for all $i$, counting multiplicities. Then $C_p(m,N,r)\!=\!0$ (if $m_i\!\leq\!1$
 for some $i$) and\\
\mbox{}\hfill\scalebox{1}[1]{$\displaystyle \pmb{C_p(m,N,r)\!\leq\!
\left\lfloor c^n
\prod\limits^n_{i=1} \left\{\left.
(m_i-1)\left[\left(\sum\limits_{j\in N_i} r_j\right) +
\log_p\left(\frac{(m_i-1)^{\#N_i}}{\left(\prod_{j\in N_i}r_j
\right)\log^{\#N_i} p} \right)\right]\right/r_i\right\}\right\rfloor}$}\hfill
\mbox{}\\ (if $m_1,\ldots,m_n\!\geq\!2$), where $\#$ denotes the operation
of taking set cardinality. \qed
\end{pcfew}

In essence, our approach replaces the use of Arithmetic Multivariate 
Descartes' Rule in \cite{amd} by simple and direct combinatorial 
argument. Furthermore, we need the $p$-adic Complex Fewnomial Theorem  
only for Theorem \ref{thm:cx}.

\section{Background and Useful Tools} 
\label{sec:back} 
Let $p$ be any prime number in the ring of ordinary integers $\Z$ and 
recall that $\Z_p$, the ring of {\bf $\pmb{p}$-adic integers}, 
can be identified with the set of all ``leftwardly infinite'' sequences of 
the form $\cdots d_2 d_1 d_0$, where $d_i\!\in\!\{0,\ldots,p-1\}$ for 
all $i\!\geq\!0$.  
(In particular, $\Z$ embeds naturally in $\Z_p$ as the set of 
all $p$-adic integers having a {\bf finite} $p$-adic expansion.) 
Addition and multiplication are then performed just as with base $p$ integers, 
noting that the carries may propagate infinitely. Put another way, 
the last $n$ digits of any $p$-adic integer calculation can be determined 
simply by working in the ring $\Z/p^{n+1}\Z$.

The $p$-adic rational numbers are then obtained by allowing 
finitely many (base $p$) digits ``after the decimal point''. 
For instance, the ordinary rational number $\frac{12345}{49}$ can  
be considered a $7$-adic rational by expressing it as the sequence of 
$7$-adic digits $506.64$, which in turn should be interpreted as 
$ 5\cdot 7^2 +6 +\frac{6}{7}+ \frac{4}{49}$. 

Following the notation of Definition \ref{dfn:padic} of the last section, 
it is then easy to show that $|x-y|_p$ defines a metric on 
$\Q_p$ and that $\Q_p$ is {\bf complete} with respect to this 
metric.\footnote{This means that all Cauchy sequences converge, 
i.e., if $\{a_i\}$ satisfies $\max\limits_{i,j\geq N} \{|a_i-a_j|_p\}
\longrightarrow 0$ as $N\longrightarrow \infty$ then 
$\{a_i\}$ has a well-defined limit. \cite{koblitz}} So the fields $\Q_p$ 
form an alternative collection of metric completions of 
$\Q$; the real numbers $\R$ being the completion of $\Q$ 
we usually see first in school. 

It is also convenient to define the $p$-adic analogue of the 
complex numbers: We let $\C_p$ denote the {\bf completion}, with 
respect to the obvious extensions of $|\cdot |_p$, of the 
algebraic closure of $\Q_p$. It turns out that 
$\C_p$ is itself algebraically closed, so no further closures 
or completions need be taken. In particular, $\ord_p \C_p\!=\!\Q$ 
\cite{koblitz}. 

One of the most amusing facts about the $p$-adic numbers is that 
their behavior as roots of polynomials can be described quite elegantly  
in terms of polyhedral geometry. 
\begin{dfn} 
Given any polynomial $f(x_1)\!=\!\sum\limits_{a\in A} c_a x^a_1
\!\in\!\C_p[x_1]$, its 
$p$-adic Newton polygon, $\newt_p(f)$, is the convex hull of\footnote{i.e., 
smallest convex set containing...} $\{(a,\ord_p c_a)\; | \; a\!\in\!A\}$ 
in $\R^2$. Also, given any polygon $Q\!\in\!\R^2$ and 
$w\!=\!(w_1,w_2)\!\in\!\R^2$, we let its 
{\bf face with inner normal $w$}, $Q^w$, be the set 
$\{(a_1,a_2)\!\in\!Q \; | \; a_1w_1+a_2w_2 \text{ minimal}\}$. 
In particular, we call $Q^w$ a vertex (resp.\ improper, an edge) 
iff $\Q^w$ is a point (resp.\ $Q$, not $Q$ or a point). 
Finally, the {\bf lower hull} of $Q$ is the union of all 
faces of $Q$ with an inner normal of the form $(v,1)$. \dia 
\end{dfn} 
\begin{prop} 
\label{prop:newt} 
(See, e.g., \cite{koblitz}.)  
Following the notation above, the number of roots $\zeta\!\in\!\C_p$ of 
$f$ with $\ord_p \zeta\!=\!v$ is exactly the length of the 
orthogonal projection of $\newt_p(f)^{(v,1)}$ onto the 
$x$-axis. \qed 
\end{prop} 

By exploiting the definition of a vertex one can easily prove the following 
useful facts on Newton polygons. 
\begin{dfn}
\label{dfn:fan} 
For any polygon $Q\!\subset\!\R^2$, let $E_+(Q)$ --- the 
{\bf unit lower edge normal locus} --- denote 
the set of $w$ in the unit circle that are inner normals 
of a lower edge of $Q$. \dia 
\end{dfn} 
\begin{prop}
\label{prop:easy}
Given any prime $p$ and $f,g\!\in\!\C_p[x_1]$, we have 
$\newt_p(fg)\!=\!\newt_p(f)+\newt_p(g)$ and 
$\newt_p(f+g)\!\subseteq\!\conv(\newt_p(f)\cup \newt_p(g))$. 
In particular, any vertex of $\newt_p(f+g)$ must be a vertex of 
$\conv(\newt_p(f)\cup \newt_p(g))$ and  
$E_+(\newt(fg))\!=\!E_+(\newt(f))\cup E_+(\newt(g))$. \qed 
\end{prop} 

We are now ready to prove our main theorems. 

\section{The Proof of Theorem \ref{thm:padic}} 
\label{sec:padic}  
That the truth of the $\tau$-conjecture implies $\pp_\C\!\neq\!\np_\C$ 
was proved by Shub and Smale in \cite{duke}. So we need only 
prove the last implication of our theorem. 

Toward this end, suppose that the $p$-adic Digit Conjecture were 
true for some fixed prime $p$. Note that any nonzero integral root $x$ of 
$f$ must have its first nonzero base-$p$ digit in $\{1,\ldots,p-1\}$. 
Furthermore, by Theorem \ref{thm:abs}, there are clearly 
no more than $N_p(\sigma(f))\!\leq\!N_p(\tau(f))\!\leq\!
\tau(f)(\tau(f)+1)/2$ possible locations for 
this first digit of $x$. Therefore, there are at most 
$pN_p(\tau(f))(1+\tau(f))^{c_p}\!\leq\!(1+\tau(f))^{c_p+1+\log_2 p}$ 
possibilities for $x\!\in\!\Z\!\setminus\!\{0\}$. 
Since $0$ is the only remaining possibility for a root of 
$f$, we are done. \qed 

\section{The Proof of Theorem \ref{thm:abs}} 
\label{sec:abs} 
The first fundamental lemma for our approach to additive complexity is 
the following. 
\begin{lemma} 
\label{lemma:upper} 
Following the notation of Theorem \ref{thm:abs}, 
$N_p(s)\!\leq\!s(s+1)/2$.
\end{lemma} 

\noindent
{\bf Proof:} Recalling the notation of Definition \ref{dfn:add}, 
note that we can define $f(x)$ via a sequence $(X_0,\ldots,X_s)$ 
where $X_{i+1}$ is a binomial or monomial in $X_0,\ldots,X_i$ for all 
$i\!\in\{0,\ldots,s\}$, with $X_0\!=\!x$, $X_{s+1}\!=\!f(x)$, and 
$s\!=\!\sigma(f)$. In particular, 
letting $L_i$ denote the number of edges on the lower hull of 
$\newt_p(X_i)$, Proposition \ref{prop:easy} immediately yields 
a recursive inequality for $L_i$: $L_{i+1}\!\leq\!2\sum^i_{j=1} 
L_j$ for all $i\!\geq\!1$. (That $L_0\!=\!0$ and $L_1\!=\!1$ is 
easily verified directly.) This recurrence easily yields $L_s\!=\!O(3^s)$ 
and thus the same upper bound on $N_p(s)$ by Proposition \ref{prop:newt}.  
However, with a little extra effort, we can do far better. 

In particular, let $F_i\!:=\!E_+(X_i)$ for all $i\!\geq\!0$. 
Then \[F_{i+1}\!=\!E_+\left(\newt_p\left(
c_{i+1} \left(\prod\limits^{i}_{j=0}
X^{m_{j,i+1}}_j\right) + d_{i+1} \left(\prod\limits^{i}_{j=0} X^{m'_{j,i+1}}_j 
\right)\right)\right),\] which by Proposition \ref{prop:easy} is 
$E_+(\conv(P\cup Q))$ where 
$P\!:=\!\newt_p\left(c_{i+1} \left(\prod\limits^{i}_{j=0}
X^{m_{j,i+1}}_j\right)\right)$ and 
$Q\!:=\!\newt_p\left(d_{i+1} \left(\prod\limits^{i}_{j=0} X^{m'_{j,i+1}}_j
\right)\right)$. In particular, note that any edge of the 
lower hull of $\conv(P\cup Q)$ must be either (a) an edge of the lower hull of 
$P$, (b) an edge of the lower hull of $Q$, or (c) a line segment connecting 
a vertex of the lower hull of $P$ with a vertex of the lower hull of 
$Q$. Moreover, since we can order the vertices of $P$ and $Q$ in, say, 
counter-clockwise order, no more than 
$1+\#(F_0\cup \cdots \cup F_i)$ edges can occur in case (c). 

Put another way, this means that $\#(F_{i+1}\!\setminus\!
(F_0\cup \cdots \cup F_i))\!\leq\!1+\#(F_0\cup \cdots \cup F_i)$. 
So by induction, we easily obtain that $L_{i+1}\!\leq\!L_i+i+1$.  
So we are done. \qed 

We are now ready to prove Theorem \ref{thm:abs}.\\
{\bf Proof of Theorem \ref{thm:abs}:} Thanks to Lemma 
\ref{lemma:upper}, we need only prove that 
$N_p(s)\!\geq\!s$. In particular, we need 
only exhibit a polynomial $f\!\in\!\Z[x_1]$ of additive complexity 
$\leq\!s$ with at least $s$ different valuations 
for its roots in $\Z$. 
Considering $(x-1)(x-p)\cdots (x-p^{s-1})$, we are done. \qed 

\begin{rem} 
Note that we could have instead worked over any ring $R\!\subseteq\!\C_p$ 
containing $\Z$ by replacing $\sigma(f)$ with $\sigma_R(f)$. 
The obvious generalization of Theorem \ref{thm:abs} then clearly 
holds. \dia 
\end{rem}  

\section{The Proof of Theorem \ref{thm:cx}}
\label{sec:cx}
First note that by the definition of additive complexity,
$x\!\in\!\C_p$ is a root
of $f \Longrightarrow
(X_0,\ldots,X_s)\!\in\!\C^{s+1}_p$ is a geometrically
isolated root of the polynomial system $G\!=\!\bO$, where the
corresponding equations are exactly
{\small
\[ c_{s+2}\prod\limits^{s+1}_{i=1}
X^{m_{i,s+2}}_i=0\]
\[X_{2} = c_{2} X^{m_{1,2}}_1 + d_{2} X^{m'_{1,2}}_1 \]
\[ \vdots \]
\[X_{s+1} = c_{s+1} \left(\prod\limits^{s}_{i=1} X^{m_{i,s+1}}_i\right) +
d_{s+1} \left(\prod \limits^{s}_{i=1} X^{m'_{i,s+1}}_i \right),\]
}

\noindent
where $s\!:=\!\sigma(f)$, $X_1\!=\!x$, $f(x)\!=\!c_{s+1}
\prod^{s+1}_{i=1}X^{m_{i,s+2}}_i$, 
and the $c_i$, $d_i$, $c^{(j)}_i$, $m_{i,j}$, and $m'_{i,j}$ are
suitable constants. This follows easily from the fact that
corresponding quotient rings $\C_p[x]/\langle f\rangle$ and
$\C_p[X_1,\ldots,X_{s+1}]/\langle G \rangle$ are isomorphic,
thus making $\C_p[x]/\langle f\rangle$ and
$\C_p[X_1,\ldots,X_{s+1}]/\langle G \rangle$ isomorphic. 

So we now need only count the geometrically isolated roots 
$(X_1,\ldots,X_{s+1})\!\in\!\C^{s+1}_p$ of $G$, 
with $|X_1-1|_p\!\leq\!\frac{1}{p^r}$,  
precisely enough to conclude. Toward this end, note that 
by a simple rescaling (and since the $p$-adic Complex Fewnomial Theorem 
is independent of the underlying coefficients) that it 
suffices to count roots of $(X_1,\ldots,X_{s+1})\!\in\!\C^{s+1}_p$ of $G$ 
with $|X_i-1|_p\!\leq\!\frac{1}{p^r}$ for all $i$. 

Note then that the first equation of \mbox{$G\!=\!\bO$} imply that at
least one $X_i$ must be $0$. 
Note also that if one of the variables $X_1,\ldots,X_{\ell+1}$ is 
$0$, then the first $1+\ell$
equations of $G$ completely determine $(X_1,\ldots,X_{\ell+1})$.
Furthermore, by virtue of the last $s-\ell$ equations of $G$,
the value of $(X_1,\ldots,X_{\ell+1})$ {\bf uniquely}
determines the value of $(X_{\ell+2},\ldots,X_{s+1})$.
So it in fact suffices to find the total number of geometrically isolated
roots (with all coordinates {\bf nonzero}) of all systems of the form 
$G'\!=\!\bO$, where the equations of $G'$ are exactly $(0\!=\!0)$ or
{\small
\[ X_{2} = c_{2} X^{m_{1,2}}_i + d_{2} X^{m'_{1,2}}_i \]
\[ \vdots \]
\[ X_{\ell} = c_{\ell} \left(\prod\limits^{\ell-1}_{i=1}
X^{m_{i,\ell}}_i\right) + d_{\ell} \left(\prod \limits^{\ell-1}_{i=1}
X^{m'_{i,\ell}}_i \right)\]
\[ 0 = c_{\ell+1} \left(\prod\limits^{\ell}_{i=1}
X^{m_{i,\ell+1}}_i\right) + d_{\ell+1} \left(\prod \limits^{\ell}_{i=1}
X^{m'_{i,\ell+1}}_i \right),\]
}

\noindent
where $\ell$ ranges over $\{1,\ldots,n\}$. Note in particular that
the $j^\thth$ equation above involves no more than $j+1$ variables for
all $j\!\in\!\{1,\ldots,\ell-1\}$, and that the $\ell^\thth$ equation
involves no more than $\ell$ variables.

Recalling the notation of the $p$-adic Complex Fewnomial Theorem, 
we then see that $G$ has no more than
\[1 \ \ , \ \ 1+C_p(2,1,r) \ \ , \ \ \rho\!:=\!1+C_p(2,1,r)+
(C_p(2,1,r)C_p(3,1,r)) \ \ , \ \ \text{or} \]
$\rho+\sum\limits^s_{\ell=3}
C_p((2,\underset{\ell-1}
{\underbrace{3,\ldots,3}}),(2,3,\ldots,\ell,\ell),\{r\}^{\ell})$  
geometrically isolated roots in $\C^{s+1}_p$ satisfying 
$|x_i-1|_p\!\leq\!\frac{1}{p^r}$ for all $i$, 
according as $s$ is $0$, $1$, $2$, or $\geq\!3$. So, by our 
earlier observations, these quantities are upper bounds for 
$C'_p(0,r)$, $C'_p(1,r)$, $C'_p(2,r)$, and $C'_p(s,r)$ respectively. 
So by an elementary calculation, we are done. \qed 

\section{The Proof of Theorem \ref{thm:add}}
\label{sec:add}
Note that any nonzero $p$-adic rational root $x$ of
$f$ must have its first nonzero digit in $\{1,\ldots,p-1\}$.
Furthermore, by Theorem \ref{thm:abs}, there are clearly
no more than $N_p(s)\!\leq\!s(s+1)/2$ possible locations for this
first digit of $x$. Therefore, by Theorem \ref{thm:cx}, there are at most
$1+pN_p(s)C'_p(s,1)$ possibilities for $x$, so we are done. 
\qed

\bibliographystyle{acm}

\begin{thebibliography}{A}

\bibitem[BCSS98]{bcss} Blum, Lenore; Cucker, Felipe; Shub, Mike; and
Smale, Steve, {\it Complexity and Real Computation,} Springer-Verlag, 1998.       
\bibitem[BC76]{bocook} Borodin, Allan and Cook, Stephen A., {\it
``On the Number of Additions to Compute Specific Polynomials,''}
SIAM J.\ Comput.\ {\bf 5} (1976), no.\ 1, pp.\ 146--157.

\bibitem[DvdD88]{vandenef} Denef, Jan and van den Dries, Lou, {\it
``$p$-adic and Real Subanalytic Sets,''} Annals of Mathematics (2) 
{\bf 128} (1988), no.\ 1, pp.\ 79--138.

\bibitem[Gri82]{grigo} Grigor'ev, Dima Yu., {\it
``Lower Bounds in the Algebraic Complexity of Computations,''}
The Theory of the Complexity of Computations, I;
Zap.\ Nauchn.\ Sem.\ Leningrad.\ Otdel.\ Mat.\ Inst.\ Steklov (LOMI) {\bf
118} (1982), pp.\ 25--82, 214.

\bibitem[GK98]{grigokarp} Grigor'ev, Dima Yu. and Karpinski, Marek, {\it
``Computing the Additive Complexity of Algebraic 
Circuits with Root Extracting,''} 
SIAM J.\ Comput., Vol.\ 27, No.\ 3, pp.\ 694--701, June 1998. 

\bibitem[Kho91]{few} Khovanski, Askold Georgievich, {\it Fewnomials,}
AMS Press, Providence, Rhode Island, 1991.

\bibitem[Kob84]{koblitz} Koblitz, Neal I., {\it $p$-adic Numbers, $p$-adic
Analysis, and Zeta-Functions,}
$2^\nd$ ed., Graduate Texts in Mathematics, 58, Springer-Verlag,
New York-Berlin, 1984. 

\bibitem[Koi96]{koiran} Koiran, Pascal, {\it ``Hilbert's 
Nullstellensatz is in the Polynomial Hierarchy,''} DIMACS Technical 
Report 96-27,July 1996. (This preprint considerably improves the published
version which appeared Journal of Complexity {\bf 12} (1996), no.\ 4,
pp.\ 273--286.)

\bibitem[Len99]{lenstra2} Lenstra, Hendrik W., Jr., {\it ``On the
Factorization of Lacunary Polynomials,''} Number Theory in Progress,
Vol.\ 1 (Zakopane-K\'oscielisko, 1997), pp.\ 277--291, de Gruyter,
Berlin, 1999.

\bibitem[Lip88]{lipshitz} Lipshitz, Leonard, {\it ``$p$-adic
Zeros of Polynomials,''} J.\ Reine Angew.\ Math.\ {\bf 390} (1988),
pp.\ 208--214.

\bibitem[dMS96]{fux} de Melo, Wellington and Svaiter, Benar Fux, 
{\it ``The Cost of Computing Integers,''} 
Proceedings of the American Mathematical Society, Vol.\ 124, No.\ 5, 
May 1996. 

\bibitem[Mor97]{gugu} de A.\ Moreira, Carlos Gustavo T., 
{\it ``On Asymptotic Estimates for Arithmetic Cost Functions,''} 
Proceedings of the American Mathematical Society, Vol.\ 125, No.\ 2, 
February 1997, pp.\ 347--353. 

\bibitem[Ris85]{risler} Risler, Jean-Jacques, {\it ``Additive Complexity
and Zeros of Real Polynomials,''} SIAM J.\ Comput.\ {\bf 14} (1985), no.\ 1,
pp.\ 178--183.

\bibitem[Roj00]{four} Rojas, J.\ Maurice, {\it ``Algebraic
Geometry Over Four Rings and the Frontier to Tractability,''} Contemporary
Mathematics, vol.\ 270, Proceedings of a Conference on Hilbert's Tenth Problem
and Related Subjects (University of Gent, November 1--5, 1999), edited by
Jan Denef, Leonard Lipschitz, Thanases Pheidas, and Jan Van
Geel, pp.\ 275--321, AMS Press (2000).

\bibitem[Roj02]{add} \underline{\hspace{\jmr}}, {\it ``Additive
Complexity and the Roots of Polynomials Over Number Fields and $\cp$-adic 
Fields,''} Proceedings of the 5$^\thth$ Annual Algorithmic Number Theory 
Symposium (ANTS V), Lecture Notes in Computer Science \#2369, pp.\ 506--515,
Springer-Verlag (2002).

\bibitem[Roj03a]{amd} \underline{\hspace{\jmr}}, {\it ``Arithmetic
Multivariate Descartes' Rule,''} Math ArXiV preprint
{\tt math.NT/0110327}, submitted for publication.

\bibitem[Roj03b]{dzh} \underline{\hspace{\jmr}}, {\it ``Dedekind Zeta 
Functions and the Complexity of Hilbert's Nullstellensatz,''} 
Math ArXiV preprint {\tt math.NT/0301111}, submitted for publication.

\bibitem[Shu93]{shub} Shub, Mike, {\it ``Some Remarks
on B\'ezout's Theorem and Complexity Theory,''} {}From
Topology to Computation: Proceedings of
the Smalefest (Berkeley, 1990), pp.\ 443--455, Springer-Verlag, 1993.

\bibitem[SS95]{duke} Shub, Mike and Smale, Steve, {\it ``On the
Intractibility of Hilbert's Nullstellensatz and an Algebraic Version of
`NP not equal to P?',''} Duke Math J., Vol.\ 81 (1995), pp.\ 47--54.

\end{thebibliography}

\end{document}